\documentclass[a4paper,11pt]{article}
\usepackage[ansinew]{inputenc}
\usepackage[english]{babel}
\usepackage{mathrsfs}
\usepackage{amsfonts}
\usepackage[centertags]{amsmath}
\usepackage{amsthm}
\usepackage{latexsym,amsmath,amssymb}
\usepackage{graphicx,color}
\usepackage{eqnarray}
\usepackage[all]{xy}% Para los diagramas.
\input diagxy
\newtheorem{teo}{Theorem}[section]
\newtheorem{prop}[teo]{Proposition}
\newtheorem{defin}[teo]{Definition}

\newtheorem{corol}[teo]{Corollary}
\newtheorem{state}[teo]{Problem}

\newtheoremstyle{obs}% name
  {3pt}%      Space above
  {3pt}%      Space below
  {}%         Body font
  {}%         Indent amount (empty = no indent, \parindent = para indent)
  {\bfseries}% Thm head font
  {.}%        Punctuation after thm head
  {.5em}%     Space after thm head: " " = normal interword space;
  {}%         Thm head spec (can be left empty, meaning `normal')
\theoremstyle{obs}

\newtheorem{remark}[teo]{Remark}

\newcommand{\ds}{\displaystyle}

\def\tabaddress#1{{\small\it\begin{tabular}[t]{c}#1
\\[1.2ex]\end{tabular}}}

\title{Strict abnormal extremals in nonholonomic and kinematic control systems}
\author{{\sc Mar\'ia Barbero-Li\~n\'an}
\thanks{{\bf e}-{\it mail}: mbarbero@ma4.upc.edu}, {\sc Miguel C. Mu\~noz-Lecanda\thanks{{\bf e}-{\it mail}:
matmcml@ma4.upc.edu}}
\\
 \tabaddress{Departamento de Matem\'atica Aplicada IV\\
  Edificio C-3, Campus Norte UPC.
  C/ Jordi Girona 1. E-08034 Barcelona. Spain}
}

\date{\today}

\begin{document}
\maketitle \pagestyle{myheadings}

\thispagestyle{empty}

\begin{abstract}
In optimal control problems, there exist different kinds of
extremals, that is, curves candidates to be solution:
abnormal, normal and strictly abnormal. The key point for
this classification is how those extremals depend on the
cost function. We focus on control systems such as
nonholonomic control mechanical systems and the associated
kinematic systems as long as they are equivalent.

With all this in mind, first we study conditions to relate
an optimal control problem for the mechanical system with
another one for the associated kinematic system. Then,
Pontryagin's Maximum Principle will be used to connect the
abnormal extremals of both optimal control problems.

An example is given to glimpse what the abnormal solutions
for kinematic systems become when they are considered as
extremals to the optimal control problem for the
corresponding nonholonomic mechanical systems.
\end{abstract}

\bigskip

\vbox{\raggedleft AMS s.\,c.\,(2000): 34A26,  49J15, 49K15,
70G45, 70H05}\null

\section{Introduction}

The problem of shortest paths in subRiemannian geometry has
strict abnormal minimizers \cite{LS96,M94}. That is why the
question of the existence of strict abnormal minimizers for
optimal control problems associated to mechanical systems
is posed here.

It will be useful to take advantage of the strict abnormal
minimizers in subRiemannian geometry to characterize, at
least, the abnormal extremals for some mechanical control
systems. In particular, we focus on the nonholomonic ones.
They are equivalent to kinematic systems under some
assumptions related to the constraint distribution and the
distribution spanned by the input vector fields, see for
instance
\cite{2005KinemBulloAndrew,2005BulloLewisBook,2007Bcn}. The
controls in those mechanical systems are understood as the
accelerations, while in the kinematic system the controls
are the velocities. The control system in subRiemannian
geometry is control-linear as the kinematic systems. From
here we connect with the nonholonomic control mechanical
system that is the object of study.

Once the equivalence between the mechanical system and the
kinematic system is established, we wonder if it is
feasible to get a similar connection between optimal
control problems associated to the two control systems. If
so, the result will be used to try to characterize strict
abnormal extremals of optimal control problems for the
nonholonomic mechanical systems.

Pontryagin's Maximum Principle
\cite{2008PMPMiguelMaria,2005BulloLewis,1997Jurdjevic,67LeeMarkus,2006AndrewCourse,P62}
defines the different kinds of extremals in optimal control
theory: normal, abnormal and strict abnormal. This
principle gives necessary conditions to find solutions to
optimal control problems. Any curve satisfying these
necessary conditions is called extremal. The extremals are
abnormal when they only depend on the geometry of the
system, in other words, the cost function does not play any
role. On the other hand, the cost function takes part in
the study of the normal extremals. In order to get a better
idea, it is said that Pontryagin's Maximum Principle
associates each solution of the optimal control problem
with a lift on the cotangent bundle, but this lift of the
solution is not necessarily unique. The non-uniqueness
makes possible the existence of extremals being normal and
abnormal at the same time. Then, a strict abnormal extremal
is one that is only abnormal, that is, it only admits one
kind of lift.

Moreover, the approach to control mechanical system
explained here enlightens how to construct the extended
system used in \cite{P62} for mechanical systems. In
contrast with the work in \cite{2005BulloLewis} our modus
operandi preserves the second order condition of the
extended system, condition also satisfied for the
non-extended system.

The paper is organized as follows: In
\S\ref{OCPnonholonkin} the different definitions and
results associated with the optimal control problems for
nonholonomic and kinematic systems are described, in
particular, the possible equivalence between both problems.
After explaining Pontryagin's Maximum Principle in
\S\ref{PMP}, the hamiltonian problems for both control
problems are stated in \S \ref{Hnonholonkin} to be able to
apply the Principle. It is especially important the
definition of extremals for the mechanical case that gives
a justification of the study made in \cite{2005BulloLewis}.
In \S \ref{example} it is showed how to use the strict
abnormal minimizers in subRiemannian geometry to
characterize the extremals for the corresponding optimal
control problem with nonholonomic mechanical system. In the
example a strict abnormal minimizer for the time optimal
control problem for the mechanical system is obtained.

 In the sequel, all the manifolds are real, second
countable and ${\cal C}^{\infty}$. The maps are assumed to
be ${\cal C}^{\infty}$. Sum over repeated indices is
understood.

\section{Optimal control problem with nonholonomic mechanical systems versus kinematic
systems}\label{OCPnonholonkin}
\subsection{Nonholonomic mechanical systems with
control}\label{S21nonholon}

Let $(Q,g)$ be a Riemannian manifold of dimension $n$ and
$\nabla$ be the Levi-Civita connection associated to the
Riemannian metric $g$. Let $TQ$ be the tangent bundle with
the natural projection $\tau_Q\colon TQ \rightarrow Q$.
Consider $D\subset TQ$, a nonintegrable distribution in $Q$
with rank $m$ and spanned by the \textit{input control
vector fields} $\{Y_1, \ldots, Y_m\}$.

Let $D^{\perp}$ be the orthogonal distribution to $D$
according to the metric $g$. Assume that $D^{\perp}$ is
spanned by $\{Z_1, \ldots, Z_{n-m}\}$, a family of vector
fields on $Q$.

It is also possible to consider an external vector field
$F\in \mathfrak{X}(Q)$, the set of smooth vector fields on
$Q$. Then, a nonholonomic mechanical system with control is
given by $\Sigma=(Q,g,F,D)$. A differentiable curve
$\gamma\colon I \rightarrow Q$ is a \textit{solution of
$\Sigma$} for certain values of the \textit{control
functions} $u^i\in \mathcal{C}^{\infty}(TQ)$ if it
satisfies the conditions
\begin{eqnarray}
\nabla_{\dot{\gamma}}\dot{\gamma}&= &F\circ \gamma
+\sum_{r=1}^{n-m}\lambda^r Z_r\circ \gamma +\sum_{s=1}^m u^s Y_s
\circ \gamma \ , \label{mechsystem}
\\\dot{\gamma}& \in & D \ . \nonumber
\end{eqnarray}
where $u\colon TQ \rightarrow U \subset \mathbb{R}^m$ being
$U$ an open set. The Lagrange multipliers $\lambda^j$ are
determined by the condition $\dot{\gamma}\in D$.

The dynamical equations of mechanical systems are second
order differential equations in the configuration manifold
$Q$, so they may be rewritten as first order differential
equations in $TQ$ using the following vector field along
the projection $\pi\colon TQ \times U \rightarrow TQ$
\begin{equation}\label{mechVF} Y=Z_g+F^V+\sum_{r=1}^{n-m}\lambda^r Z_r^V +\sum_{s=1}^m u^s Y_s^V
\end{equation}
where $Z_g$ is the geodesic spray associated with $g$ and $Y_s^V$ is
the vertical lift of $Y_s$, analogously for $F^V$ and $Z_r^V$. The
vector field $Y$ satisfies the second order condition.

On the other hand, a differentiable curve $\gamma \colon I
\rightarrow Q$ is a \textit{solution of the kinematic
system associated to (\ref{mechsystem})}, if there exist
$w^i\in \mathcal{C}^{\infty}(\mathbb{R})$, $i=1,\ldots,m$,
such that
\begin{equation}\label{kinsystem}
\dot{\gamma}(t)=\sum_{s=1}^m w^s(t) Y_s(\gamma(t)) \ ,
\end{equation}
that is, $\gamma$ is an integral curve of the vector field
$X=\sum_{s=1}^m w^sY_s$ with $w\colon \mathbb{R}\rightarrow V\subset
\mathbb{R}^m$ being $V$ an open set.

The systems (\ref{mechsystem}) and (\ref{kinsystem}) are
\textit{equivalent} if and only if every solution of
(\ref{mechsystem}) is also a solution of (\ref{kinsystem}) and in
the other way round. Notice that, in spite of the equivalence of the
systems, a solution to both systems could have different control
functions, but the curve on $Q$ is exactly the same.

\remark Here, we consider the nonholonomic control system
called \textit{fully actuated} because the constraint
distribution is exactly the distribution given by the input
control vector fields. If the distribution of the input
vector fields has rank strictly less than the rank of the
constraint distribution, then we have underactuated
systems. In this case (\ref{mechsystem}) and
(\ref{kinsystem}) are not equivalent any more, but weak
equivalent. See
\cite{2003Bloch,2005KinemBulloAndrew,2005BulloLewisBook,2007Bcn}
for more details.

\begin{teo}\label{teoequiv} \textbf{\cite{2005KinemBulloAndrew,2007Bcn}}
Every fully actuated nonholonomic control system $\Sigma$
is equivalent to the associated kinematic system.
\end{teo}

\subsection{Optimal control}\label{OC}
From a control system we define an optimal control problem adding a
cost function whose integral must be minimized over solutions of the
control system. First, we consider an optimal control problem with a
nonholonomic mechanical system. The equivalence of this system with
a kinematic system, that is, a control-linear system is known by
Theorem \ref{teoequiv}. It should be useful to find a cost function
for the kinematic system such that some connection between the
optimal solutions to both problems may be established.

Let us point out the importance of this relation between those
optimal control problems. To deal with a kinematic system is by far
easier than to deal with a mechanical control system, which is
either control-affine or nonlinear. Moreover, in \cite{LS96} the
strict abnormal minimizers have been described for the problem of
shortest-paths in subRiemannian geometry. Thus it might be expected
to characterize abnormal extremals for mechanical control systems
using the well-known abnormal minimizers in subRiemannian geometry.
The control system in subRiemannian geometry is control-linear, so
it can be understood as a kinematic system that comes from a
nonholonomic mechanical control system.

Let us consider a cost function $\mathcal{F} \colon TQ
\times U \rightarrow \mathbb{R}$ for the mechanical control
system. The optimal control problem for (\ref{mechVF}) is
stated as follows.

\begin{state}\label{mOCP} Given $x_0, x_f\in Q$, find $(\gamma,u)\colon I \rightarrow Q \times U$ such
that
\begin{enumerate}
\item $\gamma$ satisfies the end-point conditions on $Q$, i.e.
$\gamma(t_0)=x_0$, $\gamma(t_f)=x_f$;
\item $\dot{\gamma}$ is an integral curve of $Y$, i.e.
$\ddot{\gamma}(t)=Y(\dot{\gamma}(t),u(t))$;
\item $(\dot{\gamma},u)$ gives the minimum of $\int_I
\mathcal{F}(\dot{\gamma}(t),u(t))dt$ among all the curves satisfying
$1$ and $2$.
\end{enumerate}
\end{state}
In optimal control theory, it is common to consider the
functional to be minimized as a new coordinate of the
system. In this way, all the elements in the optimal
control problem are included in a control system, usually
called the \textit{extended system} \cite{67LeeMarkus,P62}.
Nevertheless, the minimization of the functional must be
included to the extended system, what turns out to be the
minimization of the new coordinate.

In the case of mechanical control systems two new coordinates are
added in order to maintain the second order condition of the vector
field (\ref{mechVF}). Let $\widehat{Q}=\mathbb{R}\times Q$, then the
cost function is considered as a vector field along the projection
$\widehat{\pi}\colon T\widehat{Q}\times U \rightarrow T\widehat{Q}$
with local expression $ \mathcal{F}{\partial}/{\partial x^0}$. Then
(\ref{mechsystem}) becomes
\begin{equation*}
\widehat{\nabla}_{\dot{\widehat{\gamma}}}
\dot{\widehat{\gamma}}= F\circ \widehat{\gamma} +
\sum_{r=1}^{n-m}\lambda^r Z_r\circ \widehat{\gamma}
+\sum_{s=1}^m u^s Y_s \circ
\widehat{\gamma}+\mathcal{F}\circ
(\dot{\widehat{\gamma}},u) \left.\frac{\partial}{\partial
x^0}\right|_{\dot{\widehat{\gamma}}} \, \end{equation*}
where $\widehat{\gamma}\colon I \rightarrow \widehat{Q}$ is
a differentiable curve and the Levi-Civita connection is
extended to $\widehat{Q}$ considering all the new
Christoffel symbols equal to zero, and $\pi_2\circ
\dot{\widehat{\gamma}}=\dot{\gamma} \in D$ with the
projection $\pi_2\colon T\widehat{Q}=T\mathbb{R} \times
TQ\rightarrow TQ$.

The above second order differential equation admits a first order
differential equation given by the vector field
\begin{equation}\label{extmechVF}
\widehat{Y}=v^0\frac{\partial}{\partial x^0}
+\mathcal{F}\frac{\partial}{\partial
v^0}+Z_g+F^V+\sum_{r=1}^{n-m}\lambda^r Z_r^V +\sum_{s=1}^m
u^s Y_s^V
\end{equation}
along the projection $\widehat{\pi}\colon
T\widehat{Q}\times U \rightarrow T\widehat{Q}$. The
differential equations added to (\ref{mechVF}) are
\[\begin{array}{rcl} \dot{x}^0 &=& v^0 \\
\dot{v}^0&=& \mathcal{F}
\end{array}\]
taking into account the extension of the Levi-Civita connection to
$\widehat{Q}$. The value that must be minimized in the optimal
control problem is $v^0=\int_I \mathcal{F} dt$.

Now consider the kinematic system $(\ref{kinsystem})$ with
a cost function $\mathcal{G} \colon Q \times V \rightarrow
\mathbb{R}$ such that the problem to be solved is
\begin{state}\label{kOCP} Given $x_0, x_f\in Q$, find $(\gamma,w)\colon I \rightarrow Q \times V$ such
that
\begin{enumerate}
\item $\gamma$ satisfies the end-point conditions on $Q$, i.e.
$\gamma(t_0)=x_0$, $\gamma(t_f)=x_f$;
\item $\gamma$ is an integral curve of $X=\sum_{s=1}^m w^s
Y_s$, i.e. $\dot{\gamma}(t)=X(\gamma(t),w(t))$;
\item $(\gamma, w)$ minimizes $\int_I \mathcal{G}(\gamma(t),w(t))dt$
among all the curves satisfying 1 and 2.
\end{enumerate}
\end{state}
\remark\label{freeOCP} The problems \ref{mOCP} and
\ref{kOCP} are called {\sl fixed time optimal control
problems} because the domain of definition of the curves is
given. However, the {\sl free time optimal control
problems} may also be defined. They consist of having
another unknown given by the domain of the definition, that
must also be found.

As before, let us extend the control system to the manifold
$\widehat{Q}=\mathbb{R}\times Q$ such that we look for integral
curves of the vector field
\begin{equation}\label{extkinVF}
\widehat{X}=\mathcal{G}\frac{\partial}{\partial
x^0}+\sum_{s=1}^m w^s Y_s
\end{equation}
defined along $\pi_1 \colon \widehat{Q}\times V \rightarrow
\widehat{Q}$. The differential equation added to
(\ref{kinsystem}) is
\[\dot{x}^0=\mathcal{G}\]
and the value to be minimized is $x^0=\int_I \mathcal{G} dt$.

By Theorem \ref{teoequiv} we know that (\ref{mechsystem}) and
(\ref{kinsystem}) are equivalent. We are interested in establishing
a connection between
\[\left. \begin{array}{rcl} \dot{x}^0=v^0 \\ \dot{v}^0=\mathcal{F} \end{array} \right\}\]
that come from (\ref{extmechVF}) and
\[\dot{x}^0=\mathcal{G}\]
that comes from (\ref{extkinVF}).

In some sense, $\mathcal{G}=v^0=\int \mathcal{F}$, but this equality
must be well understood. Observe that $\mathcal{G}$ is a function on
$\widehat{Q}\times V$, meanwhile $\mathcal{F}$ is a function on
$T\widehat{Q}\times U$. Hence, some simplifications must be
considered. Before proceeding with the exact interpretation of
$\mathcal{G}=\int \mathcal{F}$, note we also have to check what
happens with the minimization conditions when $\mathcal{G}=\int
\mathcal{F}$, that is, if the curves minimizing $\int \mathcal{G}$
determine the curves minimizing $\int \mathcal{F}$ and/or in the
other way round.

\begin{prop}\label{oneway} Let $\mathcal{G}\colon I \times Q \rightarrow
\mathbb{R}$. If $(\dot{\gamma},u)$ is an optimal curve of a
nonholonomic mechanical control system with cost function
$\mathcal{F}=\partial \mathcal{G} / \partial t +
v^i\partial \mathcal{G} / \partial x^i =\widehat{{\rm d}
\mathcal{G}}\colon I \times TQ \rightarrow \mathbb{R}$,
then there exists $w\colon I \rightarrow V$ such that
$(\gamma,w)$ is an optimal curve of the kinematic system
with cost function $\mathcal{G}$.
\end{prop}
\begin{proof}
If $(\dot{\gamma},u)\colon I \rightarrow TQ \times U$ is an integral
curve of (\ref{mechVF}), then by Theorem \ref{teoequiv} there exist
$w\colon I \rightarrow V$ such that $(\gamma,w)$ is an integral
curve of (\ref{kinsystem}). Thus, it only remains to prove that the
optimality condition for $\mathcal{F}$ implies the optimality
condition for $\mathcal{G}$.

As $(\dot{\gamma},u)$ minimizes $\int \mathcal{F}$, then for any
other integral curve $(\widetilde{\gamma},\widetilde{u})$ of the
vector field (\ref{mechVF}) satisfying the end-point conditions we
have
\[\begin{array}{lcl}
\mathcal{G}(t,\gamma(t))-\mathcal{G}(a,\gamma(a))&=&\int_a^t
\widehat{d\mathcal
{G}}(s,\gamma(s))=\int_a^t\mathcal{F}(s,\dot{\gamma}(s)){\rm
d}s <\\ \\
&<&\int_a^t\mathcal{F}(s,\dot{\widetilde{\gamma}}(s)) {\rm
d}s=\int_a^t \widehat{d\mathcal
{G}}(s,\widetilde{\gamma}(s))=\mathcal{G}(t,\widetilde{\gamma}(t))-\mathcal{G}(a,\widetilde{\gamma}(a)).
\end{array}\]
As $\gamma$ and $\widetilde{\gamma}$ satisfy the end-point
conditions and none of the cost functions depends on the
controls, we have
\[\mathcal{G}(t,\gamma(t))<\mathcal{G}(t,\widetilde{\gamma}(t)),\] then
$\int_I\mathcal{G}(t,\gamma(t)) {\rm d}t
<\int_I\mathcal{G}(t,\widetilde{\gamma}(t)){\rm d}t$ by the
monotony property of the integral.
\end{proof}

The result just proved holds provided that the cost
function for the nonholonomic mechanical system is the
total derivative of the cost function for the kinematic
system. Observe that both cost functions are independent of
the controls.

\remark Necessary conditions for a curve to be an optimal solution
for a nonholonomic mechanical control system is to be an optimal
solution to the optimal control problem for the associated kinematic
system.

\remark The inverse implication is not necessarily true. If
$(\gamma, w)$ is an optimal curve for the kinematic system, then for
any other integral curve $(\widetilde{\gamma},\widetilde{w})$ of the
kinematic system
\[\int dt\int \mathcal{F}(t, \dot{\gamma}(t))dt=\int_I\mathcal{G}(t,\gamma(t)){\rm d}t
<\int_I\mathcal{G}(t,\widetilde{\gamma}(t)){\rm d}t=\int
dt\int \mathcal{F}(t, \dot{\widetilde{\gamma}}(t))dt \ .\]
The monotony property of the integral is satisfied only in
one direction. We should think of conditions such that
\[``\int_I f <\int_I g \;\Rightarrow \; f <g \ , \quad {\rm almost \;\;
everywhere \, (a.e.)}"\] In general, we cannot expect
better results than a.e., hence we will have optimal curves
in a weak sense. For instance, if $f$ and $g$ are both
positive or both negative, the implication is satisfied.
Moreover, if $f$ and $g$ are continuous functions, then the
inequality is satisfied everywhere.

\begin{prop}\label{timeequiv} The time optimal control problem for a nonholonomic
mechanical control system is equivalent to the optimal
control problem for the associated kinematic systems with
$\mathcal{G}=t$.
\end{prop}
\begin{proof}
The direct implication is already proved in Proposition
\ref{oneway}. Let us prove now that the optimal curves for
kinematic systems with $\mathcal{G}=t$ are optimal curves
for the time optimal problem with nonholonomic mechanical
control systems.

If $(\gamma,w)$ is a minimizer of $\int t dt=t^2/2$ satisfying the
kinematic system, then by Theorem \ref{teoequiv} there exist
$u\colon I \rightarrow U$ such that $(\dot{\gamma},u)$ is an
integral curve of the nonholonomic mechanical control system. For
any other integral curve of the kinematic system with the same
end-point conditions as $\gamma$,
\[t^2/2<\widetilde{t}^2/2.\]
As $t$, $\widetilde{t}$ are positive numbers, $t<\widetilde{t}$.
That is $(\dot{\gamma},u)$ is a minimizer of the time optimal
control problem of the statement because of the equivalence of
integral curves of (\ref{mechVF}) and (\ref{kinsystem}) given by
Theorem \ref{teoequiv} and because of the nature of the cost
function. The cost function $\mathcal{G}$ is positive, so we are in
one of the cases where the reverse implication of the monotony
property of the integral is satisfied.
\end{proof}

The optimal control problems considered in Proposition
\ref{timeequiv} are free time.

\remark \label{timetime} Indeed, it is feasible to consider
the time-optimal problem for both control systems and they
will be equivalent because the time is positive. Thus, to
minimize the time or to minimize the time square is exactly
the same. Moreover, the curve on the configuration
manifolds are related to the same curve on $Q$ since the
equations defined by (\ref{mechVF}) and (\ref{kinsystem})
also appear in the extended systems (\ref{extmechVF}) and
(\ref{extkinVF}), respectively.

The following corollary links with the fact that some
optimal control problems can be understood as time optimal
control problems, as for instance happens in the problem of
shortest paths in subRiemannian geometry \cite{LS96}.

\begin{corol}  For a nonholonomic mechanical
control system, an optimal control problem equivalent to a
time optimal control problem admits an equivalent time
optimal control problem for the associated kinematic
system.
\end{corol}

The proof of this corollary is obtained from Proposition
\ref{timeequiv} and Remark \ref{timetime}.

\section{Pontryagin's Maximum Principle}\label{PMP}

Pontryagin's Maximum Principle has been widely discussed
and used in Optimal Control Theory since the second half of
the 20th century
\cite{2005BulloLewis,1997Jurdjevic,67LeeMarkus,2006AndrewCourse,P62}.

Let $Q$ be a smooth $n$-dimensional manifold and $U\subset
\mathbb{R}^m$ a bounded subset. Let $X$ be a vector field
along the projection $\pi\colon Q \times U \rightarrow Q$.
If $(x^i)$ are local coordinates on $Q$, the local
expression of the vector field is $X=f^i
{\partial}/{\partial x^i}$ where $f^i$ are functions
defined on an open set of $Q\times U$. Given ${\mathcal F}
\colon Q\times U \rightarrow \mathbb{R}$, consider the
functional
$${\cal S}[\gamma,u]=\int_I {\mathcal F}(\gamma,u)\, dt$$
defined on curves $(\gamma, u)$ with a compact interval as
domain.

To be able to state the Maximum Principle we need a
hamiltonian formalism. Now, we define the equivalent
extended optimal control problem on
$\widehat{Q}=\mathbb{R}\times Q$ with the projection
$\widehat{\pi}\colon \widehat{Q} \times U \rightarrow
\widehat{Q}$.

Let $\widehat{X}$ be the vector field along the projection
$\widehat{\pi}\colon \widehat{Q} \times U \rightarrow
\widehat{Q}$ given by:
$$\widehat{X}(x^0,x,u)={\mathcal F}(x,u) {\partial}/{\partial
x^0}|_{(x^0,x,u)} + X(x,u),$$ where $x^0$ is the natural coordinate
on $\mathbb{R}$.

\begin{state}\textbf{(Extended Optimal Control Problem, $\mathbf{\widehat{OCP}}$)}
Given $Q$, $U$, $X$, $\mathcal{F}$, $I$, $x_0$, $x_f$. Find
$(\widehat{\gamma},u)\colon I \rightarrow \widehat{Q}\times U$ such
that
\begin{enumerate}
\item $\widehat{\gamma}$ satisfies the end-point conditions: $\widehat{\gamma}(t_0)=(0,x_0)$, $\gamma(t_f)=x_f$;
\item $\dot{\widehat{\gamma}}(t)=\widehat{X}(\widehat{\gamma}(t),u(t))$ almost everywhere $t\in
I$;
\item $\gamma^0(t_f)$ is minimum over all curves
satisfying 1 and 2.
\end{enumerate}
\end{state}

The key point for considering the extended optimal control
problem is that the functional to be minimized is the
coordinate $x^0$ in $\mathbb{R}$. This is really useful in
the proof of Pontryagin's Maximum Principle and in a first
characterization of the abnormal extremals since the
direction of decreasing of the functional is easily
identified.

From $\widehat{OCP}$, we state a hamiltonian problem that
will lead to Pontryagin's Maximum Principle.

Let $T^{\ast} \widehat{Q}$ be the cotangent bundle with its
natural symplectic structure denoted by $\omega$. For each
$u \in U$, $H^u \colon T^* \widehat{Q} \rightarrow
\mathbb{R}$ is the hamiltonian function defined by
$$H^u(\widehat{x},\widehat{p})=H(\widehat{x},\widehat{p},u)=
\langle \widehat{p},\widehat{X}(\widehat{x},u) \rangle =p_0
{\mathcal F}(x,u) +\sum_{i=1}^m p_i f^i(x,u).$$ The tuple
$(T^* \widehat{Q}, \omega, H^u)$ is a hamiltonian system.
The hamiltonian vector field associated with $H$ is a
vector field along the projection $\widehat{\pi}_1 \colon
T^*\widehat{Q} \times U \rightarrow T^* \widehat{Q}$ given
by $\widehat{X}^{T^*}$, the cotangent lift of $\widehat{X}$
\cite{2005BulloLewis}.

\begin{state}\textbf{(Hamiltonian Problem, $HP$)}
 \\
 \indent Given $\widehat{OCP}$, find $(\widehat{\sigma}, u)\colon I \rightarrow T^*\widehat{Q} \times U $ such that
\begin{enumerate}
\item if
$\widehat{\gamma}=\pi_{\widehat{Q}} \circ
\widehat{\sigma}$, $\gamma=\widehat{\pi}_2 \circ
\widehat{\gamma}$ where $\widehat{\pi}_2\colon
\widehat{Q}\rightarrow Q$, then $\widehat{\gamma}(t_0)=
(0,x_0)$ and $\gamma(t_f)= x_f$;
\item
$\dot{\widehat{\sigma}}(t)=\widehat{X}^{T^*}(\widehat{\sigma}(t),u(t))$
almost everywhere $t \in I$.
\end{enumerate}
\end{state}
Locally the curve $(\widehat{\sigma},u)$ satisfies Hamilton's
equations of the system
\newline $(T^*\widehat{Q}, \omega, H^u)$,
\[\begin{array}{lcllcl} \dot{x}^0 &=& \ds{\frac{\partial
H^u}{\partial p_0}= {\mathcal F} \quad}&
\dot{p}_0&=&\ds{-\frac{\partial H^u}{\partial x^0}= 0}
\\ &&&&& \\ \dot{x}^i &=& \ds{\frac{\partial H^u}{\partial p_i}= f^i\quad }
&\dot{p}_i&=&\ds{-\frac{\partial H^u}{\partial x^i}=-p_0
\frac{\partial {\mathcal F}}{\partial x^i}- p_j
\frac{\partial f^j}{\partial x^i}}.
\end{array}\]

Note that there is no initial condition for
$\widehat{p}=(p_0, p_1, \ldots, p_n)$ in $HP$, hence it is
not a Cauchy initial value problem. This initial condition
is chosen so that the necessary conditions of Pontryagin's
Maximum Principle are satisfied, in fact, the proof of
Theorem \ref{PMfreevariable} consists of finding a suitable
initial condition
\cite{2008PMPMiguelMaria,67LeeMarkus,2006AndrewCourse,P62}.

\begin{teo} \textbf{(Pontryagin's Maximum
Principle, PMP)}\label{PMfreevariable}
\\
Let $(\widehat{\gamma},u)\colon I \rightarrow \widehat{Q}
\times U$ be a solution of the extended optimal control
problem. Then there exists $(\widehat{\sigma}, u)\colon I
\rightarrow T^*\widehat{Q} \times U$, with fiber momenta
coordinates $\widehat{\lambda}(t) \in
T^*_{\widehat{\gamma}(t)} Q$ such that:
\begin{enumerate}
\item  $(\widehat{\sigma},
u)$ is a solution of the Hamiltonian Problem;
\item $\widehat{\gamma}=\pi_{\widehat{Q}} \circ \widehat{\sigma}$;
\item \begin{itemize}
\item[(a)] $H(\widehat{\sigma}(t),u(t))= \max_{\widetilde{u} \in U} H(\widehat{\sigma}(t),
\widetilde{u})$ almost everywhere;
\item[(b)] $\max_{\widetilde{u} \in U} H(\widehat{\sigma}(t),
\widetilde{u})= {\rm constant}$ everywhere;
\item[(c)] $(\lambda_0,\lambda(t))\neq 0$ for each  $t \in
I$.
\end{itemize}
\end{enumerate}
\end{teo}

If the domain of definition of the curves is not given,
that is, free optimal control problems, see Remark
\ref{freeOCP}, then Pontryagin's Maximum Principle provides
us the same necessary conditions as Theorem
\ref{PMfreevariable}, but it also guarantees that the
maximum of the Hamiltonian is zero everywhere.

\remark \label{Hequal0} As a consequence of conditions
$(3.a)$ and $(3.b)$ the Hamiltonian along the optimal curve
with its corresponding momenta is constant almost
everywhere $t\in I$, and in particular it is zero in free
time optimal control problems. This will be used in \S
\ref{example}.

As we said previously, the proof of Theorem
\ref{PMfreevariable} consists of choosing the initial
condition for the fibers of the cotangent bundle in a
suitable way. In fact, it is chosen such that
\begin{equation}\label{separate}
\begin{array}{rcl}\langle \widehat{\sigma}(t_f),\widehat{v}(t_f)\rangle & \leq & 0 \\
\langle \widehat{\sigma}(t_f), (-1,\textbf{0})\rangle &
\geq & 0\end{array}
\end{equation} where $\widehat{v}(t_f)$ are the
perturbation vectors given by
\begin{equation}\label{vpert}
\widehat{v}(t_f)=\widehat{X}(\widehat{\gamma}(t_f),u_{t_f})-\widehat{X}(\widehat{\gamma}(t_f),u(t_f))
\end{equation}
obtained from a determined variation of the control with
value $u_{t_f}\in U$, see
\cite{2008PMPMiguelMaria,67LeeMarkus,2006AndrewCourse,P62},
and $(-1,\textbf{0})$ is the direction of decreasing in the
functional. Both vectors are in
$T_{\widehat{\gamma}(t_f)}\widehat{Q}$. Note that the
initial condition is, indeed, final since it is taken at
final time.

Observe that Maximum Principle guarantees the existence of
a covector along the optimal curve, but it does not say
anything about the uniqueness of the covector. Actually,
this covector may not be unique. Depending on the covector
we associate with the optimal curves, different kinds of
curves are defined.
\begin{defin}  \label{definextremal}
A curve $(\widehat{\gamma}, u)\colon I \rightarrow
\widehat{Q} \times U$ for $\widehat{OCP}$ is
\begin{enumerate}
\item an \textbf{extremal} if there exist $\widehat{\sigma}\colon I \rightarrow
T^*\widehat{Q}$ such that $\widehat{\gamma}=\pi_{\widehat{Q}} \circ
\widehat{\sigma}$ and $(\widehat{\sigma}, u)$ satisfies the
necessary conditions of PMP;
\item a \textbf{normal extremal}
if it is an extremal with $\lambda_0=-1$;
\item an \textbf{abnormal extremal} if it is an extremal with
$\lambda_0=0$;
\item a \textbf{strictly abnormal
extremal} if it is not a normal extremal, but it is
abnormal.
\end{enumerate}
\end{defin}
For the abnormal extremals, $\lambda_0=0$, the cost
function disappears from the hamiltonian function. Then, it
is said that the abnormal extremals only depend on the
geometry of the control system. In contrast with the normal
and strict abnormal extremals where the cost function plays
a role. In the case of strict abnormality, the cost
function is used to prove that the extremal is not normal.

\section{Hamiltonian problems for nonholonomic mechanical
systems versus kinematic systems}\label{Hnonholonkin} In order to
make profit of the optimal control problems defined in \S\ref{OC},
let us associate them with a hamiltonian problem in the sense of
Pontryagin's Maximum Principle given in \S\ref{PMP}.

For the extended mechanical system $\widehat{Y}$ given in
(\ref{extmechVF}) we have the hamiltonian function
$H_m\colon T^*T\widehat{Q}\times U\rightarrow \mathbb{R}$
defined by
\begin{equation*}
(\widehat{\Lambda},u)\longmapsto \langle \widehat{\Lambda},
v^0\frac{\partial}{\partial
x^0}+\mathcal{F}\frac{\partial}{\partial v^0}
+Z_g+F^V+\sum_{r=1}^{n-m}\lambda^r Z_r^V +\sum_{s=1}^m u^s
Y_s^V\rangle \ .\end{equation*} The Lagrange multipliers
$\lambda^j$ are fixed because they are chosen in such a way
that the part of the geodesic spray that is not in the
distribution $D$ is deleted. Another way to consider the
Lagrange multipliers is modifying the connection, see
\cite{98Lewis}.

For simplicity, we consider the system with null connection
and without external forces. Then the Lagrange multipliers
are zero and the local expression of the hamiltonian
function is
\[H_m=p_0v^0+q_0 \mathcal{F}+p_iv^i+\sum_{s=1}^m q_iu^sY^i_s \ , \]
with Hamilton's equations
\begin{equation}\label{Heqm}\begin{array}{rclrcl} \dot{x}^0&=&v^0 &
\dot{p}_0&=&0 \\
\dot{x}^i&=&v^i &
\dot{p}_i&=&\ds{-q_0\frac{\partial \mathcal{F}}{\partial x^i}-q_ju^s\frac{\partial Y^j_s}{\partial x^i}} \\
\dot{v}^0&=&\mathcal{F} & \dot{q}_0&=&-p_0\\\dot{v}^i&=&u^sY^i_s&
\dot{q}_i&=&-p_i
\end{array}\end{equation}
where $p_i$ are the momenta of the states and $q_i$ are the
corresponding momenta to the velocities.

Observe that the Hamiltonian is autonomous. Pontryagin's
Maximum Principle for this problem tells us that the
elementary perturbation vector at time $t$ for $u_1\in U$
along an optimal curve is given by
$\widehat{Y}(\dot{\widehat{\gamma}}(t),u_1)-\widehat{Y}(\dot{\widehat{\gamma}}(t),u(t))$,
see (\ref{vpert}), \begin{equation}\label{vm}
\widehat{v}_m(t)=\sum_{i=l}^m (u_1^s-u^s(t))
Y_l^V+(\mathcal{F}(\dot{\widehat{\gamma}}(t),u_1)-\mathcal{F}(\dot{\widehat{\gamma}}(t),u(t)))\left.\frac{\partial}{\partial
v^0}\right|_{\dot{\widehat{\gamma}}(t)}. \end{equation} The
covector $\widehat{\Lambda}$ associated to the optimal
curve through Pontryagin's Maximum Principle satisfies a
separating condition analogous to (\ref{separate})
\[\begin{array}{rcl}\langle \widehat{\Lambda}(t),\widehat{v}_m(t)\rangle=\langle \widehat{q}(t),\widehat{v}_m(t)\rangle & \leq & 0 \\
\langle \widehat{\Lambda}(t),
(0,\textbf{0},-1,\textbf{0})\rangle=-q_0(t) & \geq &
0.\end{array}\] The vectors $\widehat{v}_m(t)$ and
$(0,\textbf{0},-1,\textbf{0})$ are in
$T_{\dot{\widehat{\gamma}}(t)}T\widehat{Q}$. Here we do not
use the vector $(-1,\textbf{0})$, but
$(0,\textbf{0},-1,\textbf{0})$, the direction of decreasing
in the functional $\int_I{\mathcal F}$. Remember that the
value to be minimized is $v^0$.

An analogous separating condition must be satisfied for the
vector $(-1,\textbf{0},-1,\textbf{0})$ in order not to
contradict the hypothesis of optimality in Theorem
\ref{PMfreevariable}, see
\cite{2008PMPMiguelMaria,67LeeMarkus,2006AndrewCourse,P62}
for the details of that contradiction. But if
$(-1,\textbf{0},0,\textbf{0})$, the direction of decreasing
in $x^0$, is in the same half-space as the perturbation
vectors, we do not necessarily arrive at a contradiction
because, in general, a decreasing in $x^0$ does not imply a
decreasing in $v^0$.

Thus in the mechanical case the momenta must separate all
the perturbation vectors from the vectors
$(0,\textbf{0},-1,\textbf{0})$ and
$(-1,\textbf{0},-1,\textbf{0})$, what implies the
nonpositiveness of $q_0$. Taking into account Hamilton's
equations (\ref{Heqm}), $p_0$ is constant and normalizing
can be consider to be $0$, $-1$ or $1$, then the different
possibilities for the momenta are:
\begin{enumerate}
\item $p_0=0$ and $q_0=0$. Here the cost function does not
take part in the computations. Note that in this case
$(-1,\textbf{0},0,\textbf{0})$ is in the separating
hyperplane defined by the kernel of the momenta.
\item $p_0=0$ and $q_0=-1$. Then the cost function appears
in the computations. As in previous item,
$(-1,\textbf{0},0,\textbf{0})$ is in the separating
hyperplane.
\item $p_0=-1$ and $q_0=t+A$. The separating conditions
will be satisfied depending on the value of the final time
and the constant $A$. It is necessary that $A<0$ and
$t_f\leq -A$. In this case, $(-1,\textbf{0},0,\textbf{0})$
is also separated from the perturbation vectors.
\item $p_0=1$ and $q_0=-t+A$. The separating conditions
will be satisfied depending on the value of the final time
and the constant $A$. It is necessary that $A>0$ and
$t_0\geq A+1$. In this case, $(-1,\textbf{0},0,\textbf{0})$
is contained in the half-space where the perturbation
vectors are. Thus it could be associated with a
perturbation vector, depending on the directions that are
covered by the perturbations of the controls.
\end{enumerate}
To sum up, the last two previous cases cause more
difficulty to chose the initial condition for the momenta
and the final time if a free optimal control problem is
being considered. Pontryagin's Maximum Principle guarantees
the existence of a momenta, but without determining it.
Hence, we can chose the momenta that appear in the cases 1
and 2. Under that restriction, $q_0$ is a constant that
plays a similar role that the constant in Definition
\ref{definextremal}. Moreover, our mechanical Hamiltonian
turns out to be the Hamiltonian considered in
\cite{2005BulloLewis} to apply Pontryagin's Maximum
Principle for affine connection control systems. Thus, the
framework described here guarantees that the second order
condition is satisfied in the approach given in
\cite{2005BulloLewis} because it corresponds with our case
$p_0=0$.

In the extended problem for the mechanical system we have
added two new coordinates, thus two new covectors have
appeared. If we look at Definition \ref{definextremal}, it
is not clear how to define the extremals in this case. What
we have to remember is that the abnormal extremals are
characterized only using the geometry of the system before
extending it, that is, the cost function does not play any
role in the computation of abnormal extremals. For the
mechanical Hamiltonian $H_m$, this will happen if and only
if $p_0$ and $q_0$ vanish simultaneously. Otherwise, the
extremals are normal.

\begin{defin}
A curve $(\dot{\widehat{\gamma}}, u)\colon I \rightarrow
T\widehat{Q} \times U$ for the optimal control problem
\ref{mOCP} is
\begin{enumerate}
\item a \textbf{normal extremal}
if it is an extremal with either $p_0$ being a nonzero
constant or $q_0=-1$, in the latter $p_0=0$;
\item an \textbf{abnormal extremal} if it is an extremal with
$p_0=q_0=0$;
\end{enumerate}
\end{defin}

For the kinematic system, the hamiltonian function is
\begin{equation}\begin{array}{rccl}
H_k\colon& T^*\widehat{Q}\times V &\longrightarrow &
\mathbb{R}\\
&(\widehat{a},w)&\longmapsto& \ds{\langle \widehat{a},
\mathcal{G}\frac{\partial}{\partial x^0}+\sum_{s=1}^m w^s Y_s
\rangle}\end{array}\end{equation} with local expression
\[H_k=a_0\mathcal{G}+\sum_{l=1}^m
a_iw^sY^i_s \ , \] and Hamilton's equations are given by
\begin{equation}\label{Heqk}\begin{array}{rclrcl} \dot{x}^0&=&\mathcal{G} \quad &
\dot{a}_0&=&0 \\
\dot{x}^i&=&w^sY_s^i \quad & \dot{a}_i&=&\ds{-a_0\frac{\partial
\mathcal{G}}{\partial x^i}-a_jw^s\frac{\partial Y^j_s}{\partial
x^i}}
\end{array}\end{equation}
The elementary perturbation vector along the optimal curve
at $t$ for $w_1\in V$ is
\begin{equation}\label{vk}
\widehat{v}_k(t)=\sum_{s=1}^m (w_1^s-w^s(t)) Y_s+
(\mathcal{G}(\widehat{\gamma}(t),w_1)-\mathcal{G}(\widehat{\gamma}(t),w(t)))\left.\frac{\partial}{\partial
x^0}\right|_{\widehat{\gamma}(t)}\end{equation} according
to (\ref{vpert}).
 The covector
$\widehat{a}$ defined along the optimal curve that comes
from Pontryagin's Maximum Principle satisfies a separating
condition analogous to (\ref{separate})
\[\begin{array}{rcl}\langle \widehat{a}(t),\widehat{v}_k(t)\rangle & \leq & 0 \\
\langle \widehat{a}(t), (-1,\textbf{0})\rangle=-a_0 & \geq
& 0\end{array}\]where $\widehat{v}_k(t)$ and
$(-1,\textbf{0})$ are in
$T_{\widehat{\gamma}(t)}\widehat{Q}$. Here the definitions
of extremals is exactly the same as in Definition
\ref{definextremal} because there is only one more momentum
as happens in \S \ref{PMP}.

Thus we have two different hamiltonian problems, one
defined in $T^*T\widehat{Q}\times U$ and the other one
defined in $T^*\widehat{Q}\times V$. We wonder if there is
any way to relate the momenta of both problems that not
only satisfy Hamilton's equations, but also the necessary
conditions of Pontryagin's Maximum Principle. Using the
Tulczyjew diffeomorphism $\phi_{\widehat{Q}}$ defined in
\cite{Tulczyjew} there is a natural way to go from
$T^*T\widehat{Q}$ to $T^*\widehat{Q}$ with local
expression, \begin{equation}\label{mtok}
\begin{array}{ccccc}
T^*(T\widehat{Q})& \stackrel{\phi_{\widehat{Q}}}{\longrightarrow}
&T(T^*\widehat{Q})&\stackrel{\tau_{T^*\widehat{Q}}}{
\longrightarrow} & T^*\widehat{Q} \\
(x,v,p,q) & \longmapsto & (x,q,v,p) & \longmapsto & (x,q)
\end{array}\end{equation}
and it is also possible to go in the other way round as
follows
\begin{equation}\label{ktom}\xymatrix{ \stackrel{T^*\widehat{Q}}{(x,q)}&
\stackrel{T(T^*\widehat{Q})}{(x,q,\dot{x},\dot{q})}
\ar[r]^{\txt{\small{$\phi_{\widehat{Q}}^{- 1}$}}}&
\stackrel{T^*(T\widehat{Q})}{(x,\dot{x},\dot{q},q)}
 \\  I \ar[u]
\ar[ur]&&}\end{equation} where all the coordinates are
function of $t$ and $(x,q,\dot{x},\dot{q})$ is the
canonical lift of $(x(t),q(t))$ to the tangent bundle.

From here, we could think that knowing the momenta for the
mechanical system the covector for the kinematic system is
given by the momenta of the velocities. But this is not
true in general because the momenta for the kinematic
system we are looking for must also satisfy the other
necessary conditions of Pontryagin's Maximum Principle.
Moreover, both hamiltonian functions are not exactly the
same as shows (\ref{Heqm}) and (\ref{Heqk}).

In the sequel the cost functions considered for both
problems are either equal to 1, that is time optimal
problems, or the cost function given in Proposition
\ref{oneway}.

\begin{prop}\label{momenta} Let $(\widehat{\Lambda},\dot{\widehat{\gamma}})\colon I
\rightarrow T^*(T\widehat{Q})$ be a covector curve along an
optimal solution for the nonholonomic mechanical system,
Problem \ref{mOCP}. If there exists a $t_1\in I$ such that
$\langle \widehat{q}(t_1), \widehat{v}_k(t_1)\rangle \leq
0$ for every elementary perturbation vector of the
kinematic system, then $\widehat{q}(t_1)$ is the initial
condition for the covector to solve the Hamilton's equation
of the kinematic system, being $\widehat{\gamma}$  an
extremal for the kinematic Pontryagin's Maximum Principle.
\end{prop}
\begin{proof}
As an optimal solution to the nonholonomic mechanical
system is given, by Proposition \ref{oneway} and Remark
\ref{timetime} there exist controls such that the same
curve on $\widehat{Q}$ is an optimal solution to the
kinematic system. Thus, we can apply Pontryagin's Maximum
Principle that assures the existence of kinematic momenta.
But if for some $t_1\in I$, we have $\langle
\widehat{q}(t_1), \widehat{v}_k(t_1)\rangle \leq 0$, this
$\widehat{q}(t_1)$ determines the initial condition for the
momenta to integrate Hamilton's equations such that all the
necessary conditions of kinematic Pontryagin's Maximum
Principle are satisfied. The sign of the above inequality
remains invariant because of a property of the integral
curves of the complete lift and the cotangent lift of a
vector field on $\widehat{Q}$ \cite{2008PMPMiguelMaria}.
\end{proof}
\begin{corol} The abnormal optimal curves for nonholonomic
mechanical system with covectors satisfying the hypothesis in the
above proposition determine abnormal optimal curves for the
kinematic system.
\end{corol}
\begin{proof}
The momenta of abnormal extremals for nonholonomic
mechanical system are $p_0=q_0=0$. If the hypothesis in the
previous proposition are satisfied, then the initial
condition for the momenta of the kinematic system are
$\widehat{q}(t_0)$, that is, $a_0(t_0)=0$. As $a_0$ is
constant because of Hamilton's equations (\ref{Heqk}), the
abnormal solutions for the mechanical case determine
abnormal solutions for the kinematic case using Proposition
\ref{oneway} and \ref{momenta}.
\end{proof}
\remark There is an analogous result for the normal solutions as
long as the momentum for $p_0$ is taken to be equal to $0$, that is,
if we consider the case of normal solutions for mechanical systems
with $p_0=0$ and $q_0$ to be a nonzero negative constant.

\remark Observe that the extremals for the kinematic system are
extremals for the mechanical control system. But from the kinematic
momenta is not necessarily possible to find the mechanical momenta,
as the example in \S \ref{example} shows.

\subsection{Example}\label{example}
For instance, it can be proved that the example of strict
abnormal minimizer given in \cite{LS96} understood as a
solution to a nonholonomic control mechanical system is a
strict abnormal minimizer.

Let $Q=\mathbb{R}^3$ with local coordinates $(x,y,z)$. We
consider the distribution given by
$$D={\rm ker} \omega= {\rm ker} (x^2 dy - (1-x)dz)={\rm span}\{
\partial/\partial x, (1-x)\partial/\partial y + x^2 \partial/
\partial z \}= {\rm span}\{ X,Y \}.$$

Consider the Riemannian metric on $Q$, $g=dx\otimes dx+ \psi(x)
(dy\otimes dy+dz\otimes dz)$, where $\psi(x)=((1-x)^2+x^4)^{-1}$.
Observe that $X$ and $Y$ are a $g$-orthonormal basis of sections of
$Q$.

The hamiltonian function for the time optimal control
problem for the kinematic system associated to $D$ is
\[H_k(\widehat{a},w_1,w_2)=a_0+a_1 w_1+a_2 w_2
(1-x)+a_3x^2w_2.\]
 The
curve $(\gamma,w)\colon [0,1]\rightarrow Q\times V$,
$t\mapsto (0,t,0,0,1)$ satisfying the initial conditions
$\gamma(0)=(0,0,0)$ and $\gamma(1)=(0,1,0)$ is a local
strict abnormal minimizer for the time-optimal problem. It
is impossible to find momenta with $a_0=-1$ verifying all
the necessary conditions of Pontryagin's Maximum Principle.
Let us check it, the corresponding Hamilton equations for
abnormality and normality are
$$\begin{array}{ll} \dot{x}_0=1 & \dot{a}_0=0\\ \dot{x}=w_1 & \dot{a}_1=a_2w_2-2xw_2a_3 \\
\dot{y}=w_2(1-x) & \dot{a}_2=0 \\
\dot{z}=x^2w_2 & \dot{a}_3=0
\end{array}$$
Assume that the control set is open, then the maximization
of the Hamiltonian over the controls has as necessary
conditions that $\partial H_k/\partial w_1=a_1=0$,
$\partial H_k/
\partial w_2=a_2(1-x)+a_3x^2=0$. Along the curve $\gamma$,
we have $a_1=0$ and $a_2=0$. The abnormal momenta are
$\widehat{a}\colon [0,1]\rightarrow T^*\widehat{Q}$,
$t\mapsto (0,0,0,a_3)$ along $\widehat{\gamma}(t)$ with
$a_3$ being a nonzero constant. Observe that
$H_k(\widehat{a}(t),w_1,w_2)=0$ for all $t\in [0,1]$.

For the normal case, $a_0=-1$ and the necessary conditions
for the maximization of the Hamiltonian over the controls
are the same along $\gamma$: $a_1=0$, $a_2=0$. But then
$H_k(\widehat{a}(t),w_1,w_2)=-1\neq 0$ for all $t\in [0,1]$
contradicting a necessary condition of Pontryagin's Maximum
Principle, see Remark \ref{Hequal0}. Thus, as mentioned,
$\gamma$ is a strict abnormal extremal. The local
optimality is proved in \cite{LS96}.

According to the metric, the Christoffel symbols that do
not vanish are
\[\begin{array}{rcl}
\Gamma^1_{22}=\Gamma^1_{33}&=&\frac{1-x-2x^3}{((1-x)^2+x^4)^2}\\
\Gamma^2_{12}=-\Gamma^1_{21}&=&\frac{1-x-2x^3}{(1-x)^2+x^4}\\
\Gamma^3_{13}=-\Gamma^1_{31}&=&\frac{1-x-2x^3}{(1-x)^2+x^4}
\end{array}\]
where $1$ stands for coordinate $x$ and so on. Observe that
the connection associated to the metric does not have zero
torsion.

Having this in mind, the hamiltonian function for the
mechanical system is
\[H_m(\widehat{\Lambda},u_1,u_2)=p_0v_0+q_0+p_1v_1+p_2v_2+p_3v_3+q_1(-\Gamma^1_{22}v_2^2-\Gamma^1_{33}v_3^2+u_1)+
q_2u_2(1-x)+q_3x^2u_2.\] Hamilton's equations are
$$\begin{array}{ll} \dot{x}_0=v_0 & \dot{p}_0=0 \\
\dot{x}=v_1 & \dot{p}_1=\frac{\partial \Gamma^1_{22}}{\partial x}q_1v_2^2+
\frac{\partial \Gamma^1_{33}}{\partial x}q_1v_3^2+q_2u_2-2xu_2q_3 \\
\dot{y}=v_2 & \dot{p}_2=0 \\
\dot{z}=v_3 & \dot{p}_3=0\\
\dot{v}_0=1 & \dot{q}_0=-p_0\\
 \dot{v}_1=-\Gamma^1_{22}v_2^2-\Gamma^1_{33}v_3^2+u_1 & \dot{q}_1=-p_1 \\
\dot{v}_2=u_2(1-x) & \dot{q}_2=-p_2+2q_1\Gamma^1_{22}v_2 \\
\dot{v}_3=x^2u_2 & \dot{q}_3=-p_3+2v_3\Gamma^1_{33}q_1
\end{array}$$
The strict abnormal minimizer for the kinematic system becomes the
extremal $\dot{\widehat{\gamma}}(t)=(t,0,t,0,1,0,1,0)$ for the
mechanical system. Substituting into the first column of Hamilton's
equations along $\widehat{\gamma}$ we have $u_1=1$ and $u_2=1$.

\remark  The control are different for the equivalent
control systems, as was mentioned in \S \ref{S21nonholon}.

Necessary conditions for the maximization of the
Hamiltonian $H_m$ over the controls along the extremal are
$q_1=0$ and $q_2=0$. From the second column in Hamilton's
equations we have
\[\dot{p}_1=0, \; p_1=0, \; p_2=0, \; \dot{q_3}=-p_3\]
where $p_3$ is constant. These are valid for abnormality
and normality because of the considered cost function.

The abnormal momenta, $p_0=q_0=0$, is
$\widehat{\Lambda}(t)=(0,0,0,p_3,0,0,0,-p_3t+A)$ with $p_3$
and $A$ being constants, that cannot vanish simultaneously.
If now we evaluate the Hamiltonian,
$H_m(\widehat{\Lambda}(t),u_1,u_2)=0$. Thus, the abnormal
minimizer for the kinematic system is an abnormal extremal
in the mechanical case.

Let us try to find the normal momenta, that is, either
$q_0=-1$ or $p_0=-1$. The different cases are:
\begin{enumerate}
\item $p_0=-1$ then by Hamilton's equations $q_0(t)=t+B$ with a constant $B$;
\item $p_0=0$, then $q_0=-1$.
\end{enumerate}
Thus, either
$\widehat{\Lambda}_1(t)=(-1,0,0,p_3,t+B,0,0,-p_3t+A)$ or
$\widehat{\Lambda}_2(t)=(0,0,0,p_3,-1,0,0,-p_3t+A)$ along
$\dot{\widehat{\gamma}}$. If we evaluate the Hamiltonian
$H_m$ at these covectors,
\[\begin{array}{rcl}
H_m(\widehat{\Lambda}_1(t),u_1,u_2)&=&-1+t+B,\\
H_m(\widehat{\Lambda}_2(t),u_1,u_2)&=&-1
\end{array}\]
None of the previous values are zero almost everywhere on
$[0,1]$. Thus, the strict abnormal minimizer for the
kinematic system is not a normal extremal for the
mechanical case. Therefore, we have a strict abnormal
extremal for the nonholonomic mechanical system.

As for the elementary perturbation vectors $(\ref{vm})$,
$(\ref{vk})$ along the extremals considered, we have
\begin{eqnarray*}
\widehat{v}_k(t)&=&\tilde{w}_1\frac{\partial}{\partial
x}+(\tilde{w}_2-1)\frac{\partial}{\partial
y} \\
\widehat{v}_m(t)&=&(\tilde{u}_1-1)\frac{\partial}{\partial
v_1}+(\tilde{u}_2-1)\frac{\partial}{\partial v_2}.
\end{eqnarray*}
For the momenta found, the conditions $(\ref{separate})$
are
\begin{eqnarray*}
\langle (0,0,0,a_3) ,\widehat{v}_k(t) \rangle=0&\leq 0\\
\langle (0,0,0,p_3,0,0,0,-p_3t+A) ,\widehat{v}_m(t)
\rangle=0&\leq 0.
\end{eqnarray*}
Observe that the kinematic momenta and the mechanical
momenta are related through $(\ref{mtok})$, $(\ref{ktom})$.
From a kinematic system we recover the mechanical momenta
at every time $t\in I$ when $p_3=0$ and $A=a_3$. But the
way to understand the relation is: given a time $t_1$, the
initial condition for the kinematic momenta is $-p_3t_1+A$.
After integrating Hamilton's equations, the momenta do not
necessarily satisfy the relation at every time because
mechanic and kinematic Hamilton's equations are different,
although this relation is satisfied at time $t_1$. The same
happens in the other way round from the kinematic momenta
to the mechanical momenta. Thus it is highlighted the fact
that the mapping defined using Tulczyjew's diffeomorphism
does not establish a one-to-one relation between the
momenta of both Hamilton's equations for every time.

\remark Due to Proposition \ref{timeequiv} and Remark
\ref{timetime}, the strict abnormal extremal found for the
mechanical case is also a local strict abnormal minimizer
for the time optimal control problem for the control system
given by $D$.

\section*{Acknowledgements}
We acknowledge the financial support of \emph{Ministerio de
Educaci\'on y Ciencia}, Project MTM2005-04947 and the
Network Project MTM2006-27467-E/. MBL also acknowledges the
financial support of the FPU grant AP20040096.

\end{document}